\documentclass[12pt]{amsart}

\usepackage{amsmath}
\theoremstyle{plain}
\newtheorem{Thm}{Theorem}

\begin{document}

\title[Extension of Reilly formula]
{Extension of Reilly formula with applications to eigenvalue
estimates for drifting Laplacins}

\author{Li Ma, Sheng-hua Du}

\address{Department of mathematical sciences \\
Tsinghua university \\
Beijing 100084 \\
China} \email{nuslma@gmail.com} \dedicatory{}
\date{May 26th, 2009}

\begin{abstract}

In this paper, we extend the Reilly formula for drifting Laplacian
operator and apply it to study eigenvalue estimate for drifting
Laplacian operators on compact Riemannian manifolds boundary. Our
results on eigenvalue estimates extend previous results of Reilly
and Choi and Wang.


\end{abstract}

\thanks{$^*$ The research is partially supported by the National NSF grant
of China 10631020 and SRFDP 20060003002. }
 \maketitle

\section{Introduction}
Among the important formulae in differential geometry, Reilly
formula is an important tool used to give a lower bound of
eigenvalues of Laplacian operator on a Riemannian manifold with
smooth boundary. Motivated by important work of G.Perelman
\cite{P1}, we study an extension of Reilly formula for drifting
Laplacian operator associated with weighted measure and
Bakry-Emery-Ricci tensor on a compact Riemannian manifold with
smooth boundary.

Let $(M,g)$ be a compact n-dimensional Riemannian manifold with
boundary. Let $L=\Delta $ be the Laplacian operator on the compact
Riemannian manifold $(M,g)$. Given $h$ a smooth function on $M$.
We consider the elliptic operator with drifting
$$
L_h=\Delta -\nabla h\nabla
$$
associated with the weighted volume form $dm=e^{-h}dv$. We also
call $L_h$ the h-Laplacian on $M$. Assume that
\begin{equation}\label{eq1}
-L_h u=\lambda u,
\end{equation}
with $\int u^2dm=1$ and
$$
\lambda =\int |\nabla u|^2dm.
$$

 With the help of the
Bochner formula for a smooth function $f$ (see \cite{GD} and
 \cite{SY})
 $$
\frac{1}{2}L|\nabla f|^2=|D^2f|^2+(\nabla f,\nabla
Lf)+(Ric)(\nabla f,\nabla f),
 $$
we can show the following Bochner formula for Bakry-Emery-Ricci
tensor (see
 \cite{Ma})
 \begin{equation}\label{BE}
\frac{1}{2}L_h|\nabla f|^2=|D^2f|^2+(\nabla f,\nabla
L_hf)+(Ric+D^2h)(\nabla f,\nabla f).
 \end{equation}
 We remark that the tensor $Ric^h:=Ric+D^2h$ is called Bakry-Emery-Ricci
 tensor which arise naturally from the study of Ricci solitons (\cite{H95}).

 Then we have
$$
\frac{1}{2}L_h|\nabla u|^2=|D^2u|^2-\lambda|\nabla
u|^2+(Ric+D^2h)(\nabla u,\nabla u).
$$
Recall that the second fundamental form of $\partial M$ is defined
by $h(X,Y)=g(\nabla_X\nu, Y)$ for the unit outward normal $\nu$ to
the boundary $\partial M$.

 Using the integration by part on $M$, we have the
following extension of Reilly formula.

\begin{Thm}\label{thm1} We have the following extension of Reilly
formula
\begin{eqnarray}\label{EP}
&&\int_{M}(|L_hf|^2-|D^2f|^2)dm \\
&&=\int_M Ric^h(\nabla, f\nabla
f)dm+\nonumber \\
&&\int_{\partial M}((Hf_{\nu}-\nabla h\nabla
f+\Delta_{\partial}f)f_{\nu}dm \nonumber \\
&&+ \int_{\partial
M}(h(\nabla f,\nabla
f)-<\nabla_{\partial}f,\nabla_{\partial}f_{\nu}>)dm. \nonumber
\end{eqnarray}
Here and below, the symbol $\partial$ means the operation is taken
on the induced metric on $\partial M$.
\end{Thm}

We shall apply the above result to study the eigen-value estimate
for drifting Laplacian operators on $M$. We impose either
Dirichlet boundary condition $u=0$ on $\partial M$ or the Newmann
boundary condition $\frac{\partial u}{\partial \nu}=0$, where
$\nu$ is the outer unit normal to $\partial M$. The corresponding
first eigen-value of the h-Laplacian is denoted by $\lambda_D$ or
$\lambda_N$ respectively.

\begin{Thm}\label{thm2}
 Assume that
 \begin{equation}\label{ass1}
Ric+D^2h\geq \frac{|Dh|^2}{nz}+A
 \end{equation}
 for some $A>0$ and $z>0$.

 (1) If the modified mean curvature $H-h_{\nu}$ of $\partial M$ is
non-negative, then
$$
\lambda \geq \frac{n(z+1)A}{(n(z+1)-1)}.
$$

(2) If $\partial M$ is convex, that is, the second fundamental
form (defined by $h(X,Y)=g(\nabla_X\nu, Y)$) is non-negative, then
$$
\lambda \geq \frac{n(z+1)A}{(n(z+1)-1)}.
$$
\end{Thm}

We can also prove the following result, which generalize a result
of Choi and Wang \cite{CW}.
\begin{Thm}\label{thm3}
Let $(M^n,g)$ be a closed orientable manifold with $Ric^h\geq
(n-1)K>0$. Let $h$ be a smooth function on $M$.
 Let $P\subset M$ be a embedded minimal h-hypersurface
dividing $M$ into two submanifolds $M_1$ and $M_2$ (i.e.,
$H=h_{\nu}$ where $\nu$ is the unit out normal
 to $M_1$). Then for the drifting Laplacian
 $\Delta_h:=\Delta-\nabla_Ph\nabla_P$,
$$ \lambda_1(-\Delta_P^h)\geq \frac{(n-1)K}{2}.
$$
\end{Thm}

This paper is organized as follows. In section \ref{sect2} we
prove Theorem \ref{thm1}, and Theorem \ref{thm2} is proved in
section \ref{sect3}. Theorem \ref{thm3} is proved in section
\ref{sect4}.

\section{Proof of Theorem \ref{thm1}}\label{sect2}

We prove Theorem \ref{thm1} below.
\begin{proof}

We shall integra the formula (\ref{BE}). Choose a local
orthonormal frame fields $\{e_j\}$ such that $e_n=\nu$ on the
boundary $\partial M$.
 Note that
 $$
\frac{1}{2}\int_ML_h|\nabla f|^2dm=\int_{\partial
M}f_if_{ij}\nu_jdm,
 $$
and
$$
\int_M(\nabla f,\nabla L_hf)dm=\int_{\partial
M}L_hff_j\nu_jdm-\int_M |L_hf|^2dm.
$$

Then we have
\begin{eqnarray*}
&&\int_{M}(|L_hf|^2-|D^2f|^2)dm=\int_M Ric^h(\nabla f\nabla
f)dm+\\
&&\int_{\partial M}(f_n\Delta_h f-f_if_{in})dm.
\end{eqnarray*}
Recall that
$$
\Delta_h f=\Delta f-\nabla h\nabla f.
$$
Then we have
$$
f_n\Delta_h f-f_if_{in}=-f_n\nabla h\cdot\nabla
f+\sum_{j<n}(f_{jj}f_n-f_jf_{jn}).
$$

 Now
$$
\sum_{j<n}f_{jj}=\sum_{j<n}(e_j(e_jf)-(\nabla_{e_j}e_j)f)
$$
$$
=\sum_{j<n}(\nabla^{\partial}_{e_j}e_j)f-(\nabla_{e_j}e_j)f)+\Delta_{\partial}f
$$
$$
=Hf+\Delta_{\partial}f.
$$

For $j<n$,
$$
f_{jn}=f_{nj}=e_j(e_n f)-(\nabla_{e_j}e_n)f
$$
$$
=e_j(f_n)-\sum_{k<n}h_{jk}f_k.
$$
Then  we have
$$
\sum_{j<n}f_jf_{jn}=<\nabla_{\partial}f,\nabla_{\partial}f_n>-h_{jk}f_jf_k.
$$
Putting all these together we have
\begin{eqnarray*}
&&\int_{M}(|L_hf|^2-|D^2f|^2)dm=\int_M Ric^h(\nabla f,\nabla
f)dm+\\
&&\int_{\partial M}((Hf_n-\nabla h\nabla
f+\Delta_{\partial}f)f_ndm+ \\
&&\int_{\partial M}(h(\nabla f,\nabla
f)-<\nabla_{\partial}f,\nabla_{\partial}f_{\nu}>)dm.
\end{eqnarray*}

The result follows.
\end{proof}

\section{Proof of Theorem \ref{thm2}}\label{sect3}

The idea in the proof of Theorem \ref{thm2} is similar to the one
used by Reilly in \cite{R} (see also \cite{Li}). We Use the
extension of Reilly formula to prove Theorem \ref{thm2} below.
\begin{proof}
Let $L_hu+\lambda u=\Delta_hu+\lambda u=0$. We shall integrate the
extension of Reilly formula (\ref{EP}).

 Note that
$$(a+b)^2\geq \frac{a^2}{z+1}-\frac{b^2}{z}$$ for any $z>0$. So,
we have
$$
(\Delta u)^2=(\lambda u+\nabla h\nabla u)^2\geq
\frac{\lambda^2u^2}{z+1}-\frac{|\nabla h\nabla u|^2}{z}.
$$
Then we have
\begin{equation}\label{M1}
\int_{M}(|L_hu|^2-|D^2u|^2)dm \leq \int_M
(\frac{\lambda^2zu^2}{z+1}+\frac{|\nabla f\nabla u|^2}{z})dm.
\end{equation}

Note that for both Dirichlet and Newmann Problems, we have
$$
\int_{\partial M}((Hu_n-\nabla h\nabla u+\Delta_{\partial}u)u_ndm+
\int_{\partial M}(h(\nabla u,\nabla
u)-<\nabla_{\partial}u,\nabla_{\partial}u_{\nu}>)dm
$$
$$
=\int_{\partial M}((Hu_n^2-h_nu_n^2)dm+ \int_{\partial M}h(\nabla
f,\nabla f)dm\geq 0.
$$
In the last step we have used our assumption on the geometry of
$\partial M$.

 Then by our assumption (\ref{ass1}) we have
 \begin{equation}\label{M2}
\int_M Ric^h(\nabla u,\nabla u)dm\geq
\int_M(\frac{|Dh|^2}{nz}+A)|\nabla u|^2dm.
 \end{equation}
Putting (\ref{M1}) and (\ref{M2}) together we have
$$
0\geq \frac{\lambda^2}{n(z+1)}-\lambda^2+A\lambda
$$
and
$$
\lambda \geq \frac{n(z+1)A}{(n(z+1)-1)}.
$$
The result is proved.

\end{proof}

\section{Proof of Theorem \ref{thm3}}\label{sect4}

Suppose $\Delta_P^hu+\lambda u=0$.
 Without loss of generality we
may assume that
$$
\int_{\partial M_1} h(\nabla_P u,\nabla_P u)dm\geq 0.
$$

Define $f$ on $M_1$ such that
$$
\Delta_hf=0, \ \ on \ \ M
$$
with the boundary condition $f=u$ on $\partial M_1$. By Theorem
\ref{thm1} we have
\begin{eqnarray*}
&&0\geq \int_{M_1}(-|D^2f|^2)dm=\int_M Ric^h(\nabla f,\nabla
f)dm+\\
&&\int_{\partial M_1}(Hf_n-\nabla h\nabla f+\Delta_{P}u)f_ndm\\
&&+ \int_{\partial M_1}(-<\nabla_{P}f,\nabla_{P}f_{\nu}>)dm.
\end{eqnarray*}
Note that
$$
\int_{\partial M_1}(Hf_n-\nabla h\nabla
f+\Delta_{P}u)f_ndm=\int_{\partial M_1}((H-h_n)f_n-\nabla_P
h\nabla_P f+\Delta_{P}^hu)f_ndm
$$
$$
=-\int_{\partial M_1}((\nabla_{\nu} h-H)\nabla_nf+\lambda u)f_ndm
$$
$$
=-\lambda \int_{\partial M_1}uf_ndm,
$$
and
$$
\int_{\partial
M_1}(-<\nabla_{P}f,\nabla_{P}f_{\nu}>)dm=\int_{\partial M_1}
(\Delta_P^hu, f_n)dm=-\lambda \int_{\partial M_1}uf_ndm.
$$
Compute
$$
2\int_{\partial M_1}uf_ndm=\int_{\partial
M_1}(f^2)_ndm=\int_{M_1}\Delta_h f^2dm=2\int_{M_1}|\nabla f|^2dm.
$$

Using our assumption we have
$$
0\geq ((n-1)K-2\lambda)\int_{M_1}|\nabla f|^2dm.
$$
Since $\int_{M_1}|\nabla f|^2dm>0$, we get
$$
\lambda\geq \frac{(n-1)K}{2}.
$$

\end{document}